\documentclass[10.9pt,twoside]{amsart}
\usepackage{amsmath, amsthm, amscd, amsfonts, amssymb, graphicx, color}
\usepackage[bookmarksnumbered, colorlinks, plainpages]{hyperref}

\theoremstyle{definition}

\theoremstyle{remark}

\numberwithin{equation}{section}

\begin{document}
\setcounter{page}{1}
\begin{center}
{\bf ARENS REGULARITY OF TENSOR PRODUCTS   AND WEAK AMENABILITY OF BANACH ALGEBRAS }
\end{center}

\title[]{}
\author[]{KAZEM HAGHNEJAD AZAR  }

\address{}

\dedicatory{}



\begin{abstract} In this note, we study the Arens regularity of projective tensor product $A\hat{\otimes}B$ whenever $A$ and $B$ are Arens regular. We establish some new conditions for showing that the Banach algebras $A$ and $B$ are Arens regular if and only if $A\hat{\otimes}B$ is Arens regular. We also introduce some new concepts as  left-weak$^*$-weak convergence property [$Lw^*wc-$property] and right-weak$^*$-weak convergence property [$Rw^*wc-$property]  and for Banach algebra $A$, suppose that  $A^*$ and $A^{**}$, respectively,  have $Rw^*wc-$property and $Lw^*wc-$property.  Then if $A^{**}$ is weakly amenable, it follows that  $A$ is weakly amenable. We also offer some results concerning the relation between these properties with some special derivation $D:A\rightarrow A^*$. We  obtain some conclusions in the Arens regularity  of Banach algebras.
\end{abstract} \maketitle

\section{\bf  Preliminaries and
Introduction }
\noindent Suppose that $A$ and $B$ are Banach algebras. \"{U}lger in [22], has been studied that the Arens regularity of projective tensor product $A\hat{\otimes }B$. He showed that when $A$ and $B$ are Arens regular in general, $A\hat{\otimes }B$ is not Arens regular. He introduced a new concept as biregular mapping and showed that a bounded bilinear mapping $m:A\times B\rightarrow \mathbb{C}$ is biregular if and only if   $A\hat{\otimes }B$ is Arens regular. In this paper, we establish some conditions for Banach algebras $A$ and $B$ which follows that $A\hat{\otimes }B$ is Arens regular. Conversely, we investigated if $A\hat{\otimes }B$ is Arens regular, then $A$ or $B$ are Arens regular. In section three, for Banach $A-module$ $B$, we introduce  new concepts as  $left-weak^*-weak$ convergence property [ $Lw^*wc-$property] and $right-weak^*-weak$ convergence property [ $Rw^*wc-$property] with respect to $A$ and we show that if $A^*$ and $A^{**}$, respectively,  have $Rw^*wc-$property and $Lw^*wc-$property and  $A^{**}$ is weakly amenable, then $A$ is weakly amenable.  We have also some conclusions regarding Arens regularity  of Banach algebras.
We introduce some notations and definitions that we used
throughout  this paper.\\
\noindent Let $A$ be a Banach algebra and let $B$ be a   Banach $A-bimodule$.
   A derivation from $A$ into $B$ is a bounded linear mapping $D:A\rightarrow B$ such that $$D(xy)=xD(y)+D(x)y~~for~~all~~x,~y\in A.$$
The space of all continuous derivations from $A$ into $B$ is denoted by $Z^1(A,B)$.\\
Easy example of derivations are the inner derivations, which are given for each $b\in B$ by
$$\delta_b(a)=ab-ba~~for~~all~~a\in A.$$
The space of inner derivations from $A$ into $B$ is denoted by $N^1(A,B)$.
The Banach algebra $A$ is   amenable, when for every Banach $A-bimodule$ $B$, the only inner derivation from $A$ into $B^*$ is zero derivation.  It is clear that $A$ is amenable if and only if $H^1(A,B^*)=Z^1(A,B^*)/ N^1(A,B^*)=\{0\}$. The concept of amenability for a Banach algebra $A$, introduced by Johnson in 1972, has proved to be of enormous importance in Banach algebra theory, see [13].
A Banach algebra $A$ is said to be a weakly amenable, if every derivation from $A$ into $A^*$ is inner. Equivalently, $A$ is weakly amenable if and only if $H^1(A,A^*)=Z^1(A,A^*)/ N^1(A,A^*)=\{0\}$. The concept of weak amenability was first introduced by Bade, Curtis and Dales in [2] for commutative Banach algebras, and was extended to the noncommutative case by Johnson in [14].\\
 Let $A$ be  a Banach algebra and $A^*$,
$A^{**}$, respectively, be the first and second dual of $A$.  For $a\in A$
 and $a^\prime\in A^*$, we denote by $a^\prime a$
 and $a a^\prime$ respectively, the functionals in $A^*$ defined by $\langle   a^\prime a,b\rangle=\langle   a^\prime,ab\rangle=a^\prime(ab)$ and $\langle   a a^\prime,b\rangle=\langle   a^\prime,ba\rangle=a^\prime(ba)$ for all $b\in A$.
   The Banach algebra $A$ is embedded in its second dual via the identification
 $\langle   a,a^\prime\rangle$ - $\langle   a^\prime,a\rangle$ for every $a\in
A$ and $a^\prime\in
A^*$.   We
say that a bounded net $(e_{\alpha})_{{\alpha}\in I}$ in $A$ is a left bounded
approximate identity $(=LBAI)$ [resp. right bounded
approximate identity $(=RBAI)$] if,
 for each $a\in A$,   $e_{\alpha}a\longrightarrow a$ [resp. $ae_{\alpha}\longrightarrow a$].\\
 \noindent Let $X,Y,Z$ be normed spaces and $m:X\times Y\rightarrow Z$ be a bounded bilinear mapping. Arens in [1] offers two natural extensions $m^{***}$ and $m^{t***t}$ of $m$ from $X^{**}\times Y^{**}$ into $Z^{**}$ as following\\
1. $m^*:Z^*\times X\rightarrow Y^*$,~~~~~given by~~~$\langle   m^*(z^\prime,x),y\rangle=\langle   z^\prime, m(x,y)\rangle$ ~where $x\in X$, $y\in Y$, $z^\prime\in Z^*$,\\
 2. $m^{**}:Y^{**}\times Z^{*}\rightarrow X^*$,~~given by $\langle   m^{**}(y^{\prime\prime},z^\prime),x\rangle=\langle   y^{\prime\prime},m^*(z^\prime,x)\rangle$ ~where $x\in X$, $y^{\prime\prime}\in Y^{**}$, $z^\prime\in Z^*$,\\
3. $m^{***}:X^{**}\times Y^{**}\rightarrow Z^{**}$,~ given by~ ~ ~$\langle   m^{***}(x^{\prime\prime},y^{\prime\prime}),z^\prime\rangle$  $=\langle   x^{\prime\prime},m^{**}(y^{\prime\prime},z^\prime)\rangle$\\ ~where ~$x^{\prime\prime}\in X^{**}$, $y^{\prime\prime}\in Y^{**}$, $z^\prime\in Z^*$.\\
The mapping $m^{***}$ is the unique extension of $m$ such that $x^{\prime\prime}\rightarrow m^{***}(x^{\prime\prime},y^{\prime\prime})$ from $X^{**}$ into $Z^{**}$ is $weak^*-to-weak^*$ continuous for every $y^{\prime\prime}\in Y^{**}$, but the mapping $y^{\prime\prime}\rightarrow m^{***}(x^{\prime\prime},y^{\prime\prime})$ is not in general $weak^*-to-weak^*$ continuous from $Y^{**}$ into $Z^{**}$ unless $x^{\prime\prime}\in X$. Hence the first topological center of $m$ may  be defined as following
$$Z_1(m)=\{x^{\prime\prime}\in X^{**}:~~y^{\prime\prime}\rightarrow m^{***}(x^{\prime\prime},y^{\prime\prime})~~is~~weak^*-to-weak^*-continuous\}.$$
Let now $m^t:Y\times X\rightarrow Z$ be the transpose of $m$ defined by $m^t(y,x)=m(x,y)$ for every $x\in X$ and $y\in Y$. Then $m^t$ is a continuous bilinear map from $Y\times X$ to $Z$, and so it may be extended as above to $m^{t***}:Y^{**}\times X^{**}\rightarrow Z^{**}$.
 The mapping $m^{t***t}:X^{**}\times Y^{**}\rightarrow Z^{**}$ in general is not equal to $m^{***}$, see [1], if $m^{***}=m^{t***t}$, then $m$ is called Arens regular. The mapping $y^{\prime\prime}\rightarrow m^{t***t}(x^{\prime\prime},y^{\prime\prime})$ is $weak^*-to-weak^*$ continuous for every $y^{\prime\prime}\in Y^{**}$, but the mapping $x^{\prime\prime}\rightarrow m^{t***t}(x^{\prime\prime},y^{\prime\prime})$ from $X^{**}$ into $Z^{**}$ is not in general  $weak^*-to-weak^*$ continuous for every $y^{\prime\prime}\in Y^{**}$. So we define the second topological center of $m$ as
$$Z_2(m)=\{y^{\prime\prime}\in Y^{**}:~~x^{\prime\prime}\rightarrow m^{t***t}(x^{\prime\prime},y^{\prime\prime})~~is~~weak^*-to-weak^*-continuous\}.$$
It is clear that $m$ is Arens regular if and only if $Z_1(m)=X^{**}$ or $Z_2(m)=Y^{**}$. Arens regularity of $m$ is equivalent to the following
$$\lim_i\lim_j\langle   z^\prime,m(x_i,y_j)\rangle=\lim_j\lim_i\langle   z^\prime,m(x_i,y_j)\rangle,$$
whenever both limits exist for all bounded sequences $(x_i)_i\subseteq X$ , $(y_i)_i\subseteq Y$ and $z^\prime\in Z^*$, see [5, 20].\\
 The regularity of a normed algebra $A$ is defined to be the regularity of its algebra multiplication when considered as a bilinear mapping. Let $a^{\prime\prime}$ and $b^{\prime\prime}$ be elements of $A^{**}$, the second dual of $A$. By $Goldstin^,s$ Theorem [4, P.424-425], there are nets $(a_{\alpha})_{\alpha}$ and $(b_{\beta})_{\beta}$ in $A$ such that $a^{\prime\prime}=weak^*-\lim_{\alpha}a_{\alpha}$ ~and~  $b^{\prime\prime}=weak^*-\lim_{\beta}b_{\beta}$. So it is easy to see that for all $a^\prime\in A^*$,
$$\lim_{\alpha}\lim_{\beta}\langle   a^\prime,m(a_{\alpha},b_{\beta})\rangle=\langle   a^{\prime\prime}b^{\prime\prime},a^\prime\rangle$$ and
$$\lim_{\beta}\lim_{\alpha}\langle   a^\prime,m(a_{\alpha},b_{\beta})\rangle=\langle   a^{\prime\prime}ob^{\prime\prime},a^\prime\rangle,$$
where $a^{\prime\prime}b^{\prime\prime}$ and $a^{\prime\prime}ob^{\prime\prime}$ are the first and second Arens products of $A^{**}$, respectively, see [5,  20].\\
The mapping $m$ is left strongly Arens irregular if $Z_1(m)=X$ and $m$ is right strongly Arens irregular if $Z_2(m)=Y$.\\
Regarding $A$ as a Banach $A-bimodule$, the operation $\pi:A\times A\rightarrow A$ extends to $\pi^{***}$ and $\pi^{t***t}$ defined on $A^{**}\times A^{**}$. These extensions are known, respectively, as the first (left) and the second (right) Arens products, and with each of them, the second dual space $A^{**}$ becomes a Banach algebra. In this situation, we shall also simplify our notations. So the first (left) Arens product of $a^{\prime\prime},b^{\prime\prime}\in A^{**}$ shall be simply indicated by $a^{\prime\prime}b^{\prime\prime}$ and defined by the three steps:
 $$\langle   a^\prime a,b\rangle=\langle   a^\prime ,ab\rangle,$$
  $$\langle   a^{\prime\prime} a^\prime,a\rangle=\langle   a^{\prime\prime}, a^\prime a\rangle,$$
  $$\langle   a^{\prime\prime}b^{\prime\prime},a^\prime\rangle=\langle   a^{\prime\prime},b^{\prime\prime}a^\prime\rangle.$$
 for every $a,b\in A$ and $a^\prime\in A^*$. Similarly, the second (right) Arens product of $a^{\prime\prime},b^{\prime\prime}\in A^{**}$ shall be  indicated by $a^{\prime\prime}ob^{\prime\prime}$ and defined by :
 $$\langle   a oa^\prime ,b\rangle=\langle   a^\prime ,ba\rangle,$$
  $$\langle   a^\prime oa^{\prime\prime} ,a\rangle=\langle   a^{\prime\prime},a oa^\prime \rangle,$$
  $$\langle   a^{\prime\prime}ob^{\prime\prime},a^\prime\rangle=\langle   b^{\prime\prime},a^\prime ob^{\prime\prime}\rangle.$$
  for all $a,b\in A$ and $a^\prime\in A^*$.\\

\begin{center}
\section{ \bf  Arens regularity of projective tensor product algebras  }
\end{center}

\vspace {0.4 cm}

\noindent  The tensor product, $X\otimes Y$, of the vector space $X,~Y$ can be constructed as a space of linear functional on $B(X\times Y)$, in the following way:\\
Let $x\in X$ and $y\in Y$. We denote by $x\otimes y$ the functional given by evaluation at the point $(x,y)$. In other words,
$$\langle    x\otimes y, A \rangle   =A(x,y),$$
for each bilinear from $A$ on $X\times Y$ , so the tensor product $X\otimes Y$ is the subspace of the dual of bounded bilinear forms on $X\otimes Y$, $B(X\times Y)^*$.\\
We recall that each tensor $u\in X\otimes Y$ acts as a linear functional on the space of bilinear forms and so we may define a mapping $\tilde{A}:X\otimes Y\rightarrow K$ by $u\in X\otimes Y\rightarrow \langle   A,u\rangle   \in K.$ In summary, we have
$$B(X\times Y)=(X\otimes Y)^*.$$
Let $X,Y,E~and ~F$ be vector spaces and let $S:X\rightarrow E$ and $T:Y\rightarrow F$ be linear mappings. Then we may define a bilinear mapping by $(x,y)\in X\times Y\rightarrow (Sx)\otimes (Ty)\in E \otimes F$. Linearization gives a linear mapping $(S\otimes T): X\otimes Y\rightarrow  E \otimes F$ such that             $(S\otimes T)(x\otimes y)=(Sx)\otimes (Ty)$ for every $x\in X$ and $y\in Y$.\\
By $X\hat{\otimes }Y$ and $X\check{\otimes} Y$ we shall denote, respectively, the projective and injective tensor products of $X$ and $Y$. That is, $X\hat{\otimes }Y$ is the completion of $X\otimes Y$ for the norm
$$\parallel u\parallel=inf \sum_{i=1}^n \parallel x_i\parallel\parallel y_i\parallel,$$
where the infimum is taken over all the representations of $u$ as a finite sum of the form $u=\sum_{i=1}^n x_i\otimes y_i$, and $X\check{\otimes} Y$ is the completion of $X\otimes Y$ for the norm
$$\parallel u\parallel=sup\{\mid\sum_{i=1}^n \langle   x^\prime,x_i\rangle   \langle   y^\prime,y_i\rangle   \mid:~\parallel x^\prime\parallel \leq 1,~\parallel y^\prime\parallel \leq 1\}.$$
The dual space of $X\hat{\otimes }Y$ is $B(X\times Y)$, and that of $X\check{\otimes} Y$ is a subspace of $B(X\times Y)$.\\
Although the injective tensor product of two Banach algebra $A$ and $B$ is not always a Banach algebra, their projective tensor product is always a Banach algebra. The natural multiplication of $A\hat{\otimes }B$ is the linear extension of the following multiplication on decomposable tensors $(a\otimes b)(\tilde{a}\otimes \tilde{b})=a\tilde{a}\otimes b\tilde{b}$. For more information about the tensor product of Banach algebra, see for example  [4, 5].\\
A functional $a^\prime$ in $A^*$ is said to be $wap$ (weakly almost
 periodic) on $A$ if the mapping $a\rightarrow a^\prime a$ from $A$ into
 $A^{*}$ is weakly compact. Pym in [20] showed that  this definition to the equivalent following condition\\
 For any two net $(a_{i})_{i}$ and $(b_{j})_{j}$
 in $\{a\in A:~\parallel a\parallel\leq 1\}$, we have\\
$$\\lim_{i}\\lim_{j}\langle   a^\prime,a_{i}b_{j}\rangle   =\\lim_{j}\\lim_{i}\langle   a^\prime,a_{i}b_{j}\rangle   ,$$
whenever both iterated limits exist. The collection of all $wap$
functionals on $A$ is denoted by $wap(A)$. Also we have
$a^{\prime}\in wap(A)$ if and only if $\langle   a^{\prime\prime}b^{\prime\prime},a^\prime\rangle   =\langle   a^{\prime\prime}ob^{\prime\prime},a^\prime\rangle   $ for every $a^{\prime\prime},~b^{\prime\prime} \in
A^{**}$. Thus, it is clear that $A$ is Arens regular if and only if $wap(A)=A^*$. In the following, for Banach algebras $A$ and $B$, for showing Arens regularity of projective tensor products $A\hat{\otimes}B$, we establish $wap(A\hat{\otimes}B)=(A\hat{\otimes}B)^*$. In all of this section, we regard $A^*\hat{\otimes}B^*$ as a subset of $(A\hat{\otimes}B)^*$.\\\\

\noindent{\it{\bf Theorem 2-1.}} Suppose that $A$ and $B$ are  Banach algebra and for every sequence $(x_i)_i~,(y_j)_j\subseteq A_1$, $(z_i)_i~,(w_j)_j\subseteq B_1$ and $f\in \mathbf{B}(A\times B),$ we have
$$\lim_j\lim_if(x_iz_i,y_jw_j)=\lim_i\lim_jf(x_iz_i,y_jw_j).$$
Then $A\hat{\otimes}B$ is Arens regular.\\

\begin{proof} Assume that $f\in \mathbf{B}(A\times B)$. Since $\mathbf{B}(A\times B)=(A\hat{\otimes}B)^*$, it is enough to show that $f\in wap(A\hat{\otimes}B)$.
Let $(x_i)_i~,(y_j)_j\subseteq A_1$ and  $(z_i)_i~,(w_j)_j\subseteq B_1$, then we have the following equality
$$\lim_j\lim_i\langle   f,(x_i\otimes y_j)(z_i\otimes w_j)\rangle   =\lim_j\lim_i\langle   f,x_iz_i\otimes y_jw_j\rangle   $$$$=\lim_j\lim_if(x_iz_i,y_jw_j)=\lim_i\lim_jf(x_iz_i,y_jw_j)$$$$=\lim_i\lim_j\langle   f,(x_i\otimes y_j)(z_i\otimes w_j)\rangle   .$$
Consequently by [20], $f\in wap(A\hat{\otimes}B)$.\end{proof}
\vspace {0.6 cm}

 \noindent{\it{\bf Definition 2-2.}} Assume that  $B$ is a Banach $A-bimodule$. We say that $B$ is non-trivial on $A$, if  for every  $(a_i)_{i=1}^n\subseteq A_1$ and $(b_j)_{j=1}^n\subseteq B_1$, respectively, basis elements of $A$ and $B$, we have $\sum_{i=1}^n \alpha_ia_ib_i\neq 0$  where   $\alpha_i$ is scaler and  every $a_i$ and $b_i$ are distinct for all $1\leq i\leq n$.\\
For example, take $B=\mathbb{R}\times \{0\}$ and $A=\mathbb{R}^2$ by the following multiplication
$$(a_1,a_2)(b_1,0)=(a_1b_1, 0)~~ where~~ a_1,a_2,b_1\in \mathbb{R}.$$\\

\noindent{\it{\bf Theorem 2-3.}} Suppose that $A$ and $B$ are  Banach algebras and $B$ is unital. Let $B$ be a Banach $A-bimodule$. Then we have the following assertions:
 \begin{enumerate}
 \item ~ If $A\hat{\otimes}B$ is Arens regular, then $A$ is Arens regular.
 \item ~ Let $B$ be non-trivial on $A$ and let $B$ be an unital Banach $A-module$. Then  $A$ and $B$ are Arens regular if and only if $A\hat{\otimes}B$ is Arens regular.  \end{enumerate}\begin{proof}

\vspace {0.6 cm}

 \begin{enumerate}
 \item ~ Assume that $A\hat{\otimes}B$ is Arens regular and let $u\in B$ be an unit element of $B$.  We show that $wap(A)=A^*$. Assume that   $(a_i)_i\subseteq A$ , $(c_j)_j\subseteq A$ whenever both iterated limits exist and $a^\prime\in A^*$. Then we define $\phi=a^\prime\otimes b^\prime$ where $b^\prime\in B^*$ and $b^\prime(u)=1$. Since $A^*\otimes B^*\subseteq (A\otimes B)^*$ and $A\hat{\otimes}B$ is Arens regular, we have $a^\prime\otimes b^\prime\in wap(A\otimes B)$. Hence it follows that
$$\lim_i\lim_j\langle   a^\prime, a_ic_j\rangle   =\lim_i\lim_j\langle   a^\prime\otimes b^\prime, a_ic_j\otimes u\rangle   $$$$=\lim_i\lim_j\langle   a^\prime\otimes b^\prime, (a_i\otimes u)(c_j\otimes u)\rangle   =\lim_j\lim_i\langle   a^\prime\otimes b^\prime, (a_i\otimes u)(c_j\otimes u)\rangle   $$$$=\lim_j\lim_i\langle   a^\prime, a_ic_j\rangle   .$$
We conclude that $a^\prime \in wap(A)$, and so $A$ is Arens regular.

 \item ~Let $u$ be an unit element of $B$ and suppose that   $B$ is  Arens regular. Then $wap(B)=B^*$.   Suppose that $(a_i)_i\subseteq A_1$ and $(b_j)_j\subseteq B_1$ whenever both iterated limits exist. Then $(a_iu)_i\subseteq B_1$, and so for every $b^\prime\in B^*$, we have the following equality
$$\lim_i\lim_j\langle   b^\prime,(a_iu)b_j\rangle   =\lim_j\lim_i\langle   b^\prime,(a_iu)b_j\rangle   .$$
Now let $\phi\in(A\hat{\otimes}B)^*$. We define the mapping $T:A\hat{\otimes}B\rightarrow B$ such that $T(\sum_{i=1}^n \alpha_ia_i\otimes b_i)=\sum_{i=1}^n \alpha_ia_ib_i$ where $a_i\in A$, $b_i\in B$ and $\alpha_i$ is a scaler. We show  that $\phi oT^{-1}\in B^*$. Since $B$ is not-trivial on $A$, $T^{-1}$ exist. Now let $e\in A$ be an unit element for $B$ as Banach $A-module$ and let $(b_\alpha)_\alpha\subseteq B$ such that $b_\alpha\rightarrow b$. Then $e\otimes b_\alpha\rightarrow e\otimes b$ in $A\hat{\otimes}B$, it follows that
$$\langle   \phi oT^{-1},b_\alpha\rangle   =\langle   \phi ,T^{-1}(b_\alpha)\rangle   =\langle   \phi ,T^{-1}(eb_\alpha)\rangle   =\langle   \phi ,e\otimes b_\alpha\rangle   $$$$\rightarrow
\langle   \phi ,e\otimes b\rangle   =\langle   \phi ,T^{-1}(ub)\rangle   =\langle   \phi oT^{-1},b\rangle   .$$

Consequently  $\phi oT^{-1}\in B^*$. Now we have the following equality
$$\lim_i\lim_j\langle   \phi, a_i\otimes b_j\rangle   =\lim_i\lim_j\langle   b^\prime o T, a_i\otimes b_j\rangle   $$$$=\lim_i\lim_j\langle   b^\prime, T(a_i\otimes b_j)\rangle   =
\lim_i\lim_j\langle   b^\prime, a_i b_j\rangle   $$$$=\lim_i\lim_j\langle   b^\prime,a_i(ub_j)\rangle   =\lim_j\lim_i\langle   b^\prime,(a_iu)b_j\rangle   $$$$=\lim_j\lim_i\langle   b^\prime,T(a_i\otimes b_j)\rangle   =
\lim_j\lim_i\langle   \phi, a_i\otimes b_j\rangle   .$$
It follows that $\phi\in wap(A\hat{\otimes}B)$, and so $A\hat{\otimes}B$ is Arens regular.\\
The converse by using part (1)  hold.\end{enumerate}\end{proof}

\vspace {0.6 cm}

\noindent{\it{\bf Corollary 2-4 .}} Suppose that $A$ and $B$ are unital Banach algebras and  $B$ is an unital Banach as $A-module$.    Assume that  $B$ is non-trivial on $A$. Then if  $A$ and $B$ are Arens regular, then every bilinear form $m:A \times B\rightarrow \mathbb{C}$ is weakly compact.
\begin{proof} By using Theorem 2-3 and [22, Theorem 3.4], proof  hold.\end{proof}

\vspace {0.6 cm}

\noindent{\it{\bf Example 2-5.}} $(\ell^1\oplus\mathbb{C})\hat{\otimes}\ell^\infty$ is Arens regular.
\begin{proof} We know that $\ell^\infty$ is  $(\ell^1\oplus\mathbb{C})-bimodule$ and $\ell^\infty$ is unital. $\ell^\infty$ is also non-trivial on   $(\ell^1\oplus\mathbb{C})$. By using [2, Corollary 8] and [5, Example 2.6.22(iii)], respectively, we know that $\ell^\infty$ and $(\ell^1\oplus\mathbb{C})$ are Arens regular,  and so by Theorem 2-3, $(\ell^1\oplus\mathbb{C})\hat{\otimes}\ell^\infty$ is Arens regular.\end{proof}
\vspace {0.6 cm}

 Let $A$ and $B$ be Banach algebras. A bilinear form $m:A \times B\rightarrow \mathbb{C}$ is said to be biregular,  if for any two pairs of sequence
$({a}_i)_i$ ,  $(\tilde{a}_j)_j$ in $A_1$ and $(b_i)_i$, $(\tilde {b}_j)_j$
 in $B_1$, we have
$$\lim_i\lim_jm(a_i\tilde{a}_j,b_i\tilde{b}_j)=\lim_j\lim_im(a_i\tilde{a}_j,b_i\tilde{b}_j)$$
provided that these limits exist.\\
There are some example of  biregular non regular bilinear form that for more information see [22].\\\\

 \noindent{\it{\bf Corollary 2-6.}} Suppose that $A$ and $B$ are  Banach algebras. Then we have the following assertions.
 \begin{enumerate}
 \item ~By conditions of Theorem 2-1, every bilinear form $m:A \times B\rightarrow \mathbb{C}$ is biregular.
 \item ~By conditions of Theorem 2-3 (2), every bilinear form $m:A \times B\rightarrow \mathbb{C}$ is biregular.\\\\
 \end{enumerate}

 \noindent{\it{\bf Example 2-7.}} Every bilinear form $m:(\ell^1\oplus\mathbb{C})\times\ell^\infty\rightarrow \mathbb{C}$ is Arens regular.
\begin{proof} By notice to Example 2-5 and [22, Theorem 3.4], proof is hold.\end{proof}

\vspace {0.6 cm}

 \noindent In the following we give simple proof for biregularity of bilinear form $m:A \times B\rightarrow \mathbb{C}$ such that $m(a,b)=\langle   u(a),b\rangle   $ where $u:A\rightarrow B^*$ is continuous linear operator that is introduced in [22 , Theorem 3.4].\\\\

\noindent{\it{\bf Theorem 2-8 [22].}} Let $A$ and $B$ be Banach algebras and $u:A\rightarrow B^*$ is continuous linear operator. Then the bilinear form $m:A \times B\rightarrow \mathbb{C}$ defined by  $m(a,b)=\langle   u(a),b\rangle   $ is biregular.

\begin{proof} Let $(a_i)_i$, $(\tilde{a}_j)_j$ in $A_1$ and $(b_i)_i$, $(\tilde{b}_j)_j$ in $B_1$ be such that the following iterated limits exist
$$\lim_i\lim_jm(a_i\tilde{a}_j,b_i\tilde{b}_j)~~and~~\lim_j\lim_im(a_i\tilde{a}_j,b_i\tilde{b}_j).$$
By [8, p.424], from these sequences we can extract $(a_\alpha)_\alpha$, $(\tilde{a}_\beta)_\beta$ in $A$ and $(b_\alpha)_\alpha$, $(\tilde{b}_\beta)_\beta$ in $B$ such that $a_{\alpha} \stackrel{w^*} {\rightarrow}a^{\prime\prime}$ and $\tilde{a}_{\beta} \stackrel{w^*} {\rightarrow}\tilde{a}^{\prime\prime}$ in $A^{**}$ and we have also $b_{\alpha} \stackrel{w^*} {\rightarrow}b^{\prime\prime}$ and $\tilde{b}_{\beta} \stackrel{w^*} {\rightarrow}\tilde{b}^{\prime\prime}$ in $B^{**}$. Since $A$ and $B$ are Arens regular, we have
$$\lim_\alpha \lim_\beta a_\alpha\tilde{a}_\beta=\lim_\beta \lim_\alpha a_\alpha\tilde{a}_\beta=a^{\prime\prime}\tilde{a}^{\prime\prime}~~$$ $$and$$
$$\lim_\alpha \lim_\beta b_\alpha\tilde{b}_\beta=\lim_\beta \lim_\alpha b_\alpha\tilde{b}_\beta=b^{\prime\prime}\tilde{b}^{\prime\prime}~~$$
Then, since $u$ is continuous, we have
$$\lim_\alpha \lim_\beta m(a_\alpha\tilde{a}_\beta,b_\alpha\tilde{b}_\beta)= \lim_\alpha \lim_\beta \langle   u(a_\alpha\tilde{a}_\beta),b_\alpha\tilde{b}_\beta\rangle   $$$$=\langle   u^{\prime\prime}
(a^{\prime\prime}\tilde{a}^{\prime\prime}),b^{\prime\prime}\tilde{b}^{\prime\prime}\rangle   .$$
Similarly, we have
$$\lim_\beta \lim_\alpha m(a_\alpha\tilde{a}_\beta,b_\alpha\tilde{b}_\beta)=\langle   u^{\prime\prime}
(a^{\prime\prime}\tilde{a}^{\prime\prime}),b^{\prime\prime}\tilde{b}^{\prime\prime}\rangle   .$$
Consequently we have
$$\lim_i\lim_j m(a_i\tilde{a}_j,b_i\tilde{b}_j)=\lim_j\lim_i m(a_i\tilde{a}_j,b_i\tilde{b}_j).$$
It follows that $m$ is biregular.\end{proof}

\vspace {0.6 cm}

\noindent{\it{\bf Example 2-9 [22].}} Let $A$ be a Banach algebra and $1<   p<   \infty$. Then
\begin{enumerate}
\item $\ell^p\hat{\otimes}A$ is Arens regular if and only if $A$ is Arens regular.
\item  Let $G$ be a locally compact group. Then, $L^p(G)\hat{\otimes}A$ is Arens regular if and only if $A$ is Arens regular.
\end{enumerate}
\begin{proof} By using [22, Theorem 3.4] and Theorem 2-8, proof  hold.\end{proof}

\vspace {0.6 cm}

\section{ \bf  Weak amenability of Banach algebras  }
\vspace{0.1cm}
\noindent  For Banach algebra $A$, Dales,  Rodrigues-Palacios and  Velasco in [7] have been studied the weak amenability of $A$, when its second dual is weakly amenable.   Mohamadzadih and Vishki in [19] have given simple solution to this problem with some other results, and  Eshaghi Gordji and Filali in [10] have been studied  this problem with some new results. In this section, We study this problem in the new way with some new results. Thus,   for Banach $A-module$ $B$, we introduce some new concepts as  $left-weak^*-weak$ convergence property [ $Lw^*wc-$property] and $right-weak^*-weak$ convergence property [ $Rw^*wc-$property] with respect to $A$ and we show that if $A^*$ and $A^{**}$, respectively,  have $Rw^*wc-$property and $Lw^*wc-$property and  $A^{**}$ is weakly amenable, then $A$ is weakly amenable. We also show the relations between these properties and weak amenability of $A$.
Now in the following, for left and right Banach $A-module$ $B$, we define, respectively, $Lw^*wc-$property and $Rw^*wc-$property concepts with some examples. \\

\noindent{\it{\bf Definition 3-1.}} Assume that  $B$ is a left Banach $A-module$. Let $a^{\prime\prime}\in A^{**}$ and $(a_\alpha)_\alpha \subset A$ such that $ a_\alpha  \stackrel{w^*} {\rightarrow}a^{\prime\prime}$ in $A^{**}$. We say that $b^\prime\in B^*$ has $left-weak^*-weak$ convergence property  $Lw^*wc$-property with respect to $A$, if $b^\prime a_\alpha  \stackrel{w} {\rightarrow}b^\prime a^{\prime\prime}$ in $B^*$.\\
 When every $b^\prime \in B^*$ has   $Lw^*wc$-property with respect to $A$, we say that $B^*$ has $Lw^*wc-$property.\\
 The definition of $right-weak^*-weak$ convergence property [$=Rw^*wc-$property] with respect to $A$ is similar and if $b^\prime \in B^*$ has $left-weak^*-weak$ convergence property and $right-weak^*-weak$ convergence property, then we say that $b^\prime \in B^*$ has $weak^*-weak$ convergence property  [$=w^*wc-$property].\\
 By using Lemma 3.1 from [17], it is clear that if $A^*$ has $Lw^*wc$-property, then $A$ is Arens regular.\\
 Assume that  $B$ is a left Banach $A-module$. We say that $b^\prime\in B^*$ has $left-weak^*-weak$ convergence property to zero $Lw^*wc$-property to zero with respect to $A$, if for every $(a_\alpha)_\alpha \subset A$, $b^\prime a_\alpha  \stackrel{w^*} {\rightarrow}0$ in $B^*$ implies that $b^\prime a_\alpha  \stackrel{w} {\rightarrow}0$ in $B^*$.\\

\noindent{\it{\bf Example 3-2 .}} \begin{enumerate}

 \item ~ Every reflexive Banach $A-module$ has $w^*wc-$property.

  \item ~Let $\Omega$ be a compact group and suppose that $A=C(\Omega)$ and $B=M(\Omega)$. Let $(a_\alpha)_\alpha \subseteq A$ and $\mu\in B$. Suppose that $\mu  a_\alpha  \stackrel{w^*} {\rightarrow}0$, then for each $a\in A$, we have
 $$\langle   \mu a_\alpha ,a\rangle=\langle   \mu , a_\alpha *a\rangle=\int_\Omega (a_\alpha *a) d\mu\rightarrow 0.$$
 We set $a=1_\Omega$ . Then $\mu(a_\alpha)\rightarrow 0$. Now let $b^\prime\in B^*$. Then
 $$\langle   b^\prime, \mu a_\alpha\rangle=\langle   a_\alpha b^\prime, \mu  \rangle=\int_\Omega a_\alpha b^\prime d\mu\leq \parallel b^\prime \parallel~\mid\int_\Omega a_\alpha d\mu \mid=\parallel b^\prime \parallel~\mid\mu( a_\alpha ) \mid\rightarrow 0.$$
It follows that  $\mu a_\alpha  \stackrel{w} {\rightarrow}0$, and so that $\mu$ has $Rw^*wc-$property to zero with respect to $A$.\\  \end{enumerate}

Let now $B$ be a Banach $A-bimodule$, and let\\
$$\pi_\ell:~A\times B\rightarrow B~~~and~~~\pi_r:~B\times A\rightarrow B.$$
be the left and right module actions of $A$ on $B$, respectively. Then $B^{**}$ is a Banach $A^{**}-bimodule$ with module actions
$$\pi_\ell^{***}:~A^{**}\times B^{**}\rightarrow B^{**}~~~and~~~\pi_r^{***}:~B^{**}\times A^{**}\rightarrow B^{**}.$$
Similarly, $B^{**}$ is a Banach $A^{**}-bimodule$ with module actions\\
$$\pi_\ell^{t***t}:~A^{**}\times B^{**}\rightarrow B^{**}~~~and~~~\pi_r^{t***t}:~B^{**}\times A^{**}\rightarrow B^{**}.$$\\
For a Banach $A-bimodule$ $B$,  we  define the topological centers of the  left and right module actions of $A$ on $B$ as follows:

$${Z}^\ell_{A^{**}}(B^{**})={Z}(\pi_r)=\{b^{\prime\prime}\in B^{**}:~the~map~~a^{\prime\prime}\rightarrow \pi_r^{***}(b^{\prime\prime}, a^{\prime\prime})~:~A^{**}\rightarrow B^{**}$$$$~is~~~weak^*-weak^*~continuous\}$$
$${Z}^\ell_{B^{**}}(A^{**})={Z}(\pi_\ell)=\{a^{\prime\prime}\in A^{**}:~the~map~~b^{\prime\prime}\rightarrow \pi_\ell^{***}(a^{\prime\prime}, b^{\prime\prime})~:~B^{**}\rightarrow B^{**}$$$$~is~~~weak^*-weak^*~continuous\}$$
$${Z}^r_{A^{**}}(B^{**})={Z}(\pi_\ell^t)=\{b^{\prime\prime}\in B^{**}:~the~map~~a^{\prime\prime}\rightarrow \pi_\ell^{t***}(b^{\prime\prime}, a^{\prime\prime})~:~A^{**}\rightarrow B^{**}$$$$~is~~~weak^*-weak^*~continuous\}$$
$${Z}^r_{B^{**}}(A^{**})={Z}(\pi_r^t)=\{a^{\prime\prime}\in A^{**}:~the~map~~b^{\prime\prime}\rightarrow \pi_r^{t***}(a^{\prime\prime}, b^{\prime\prime})~:~B^{**}\rightarrow B^{**}$$$$~is~~~weak^*-weak^*~continuous\}.$$
\vspace{0.4cm}

 \noindent{\it{\bf Theorem 3-3.}} i) Assume that  $B$ is a left Banach $A-module$.
If $B^*A^{**}\subseteq B^*$, then $B^*$ has $Lw^*wc-$property.\\
 ii) Assume that  $B$ is a right Banach $A-module$.
If $A^{**}B^*\subseteq B^*$  and $Z^r_{A^{**}}(B^{**})=B^{**}$, then $B^*$ has $Rw^*wc-$property.\\
\begin{proof} i) Assume that $a^{\prime\prime}\in A^{**}$ and $(a_\alpha)_\alpha \subseteq A$ such that $a_\alpha \stackrel{w^*} {\rightarrow}a^{\prime\prime}$. Then for every $b^{\prime\prime}\in B^{**}$, since $b^\prime a^{\prime\prime}\in B^*$, we have
$$<b^{\prime\prime}, b^\prime a^{\prime\prime}>=<a^{\prime\prime}b^{\prime\prime}, b^\prime >=\lim_\alpha <a_\alpha b^{\prime\prime}, b^\prime >=\lim_\alpha < b^{\prime\prime}, b^\prime a_\alpha>.$$
It follows that $b^\prime a_\alpha \stackrel{w} {\rightarrow}b^\prime a^{\prime\prime}$.\\
ii) Proof is similar to (i).\end{proof}
\vspace{0.4cm}

 \noindent{\it{\bf Theorem 3-4.}} Let $A$ be a Banach algebra and suppose that $A^*$ and $A^{**}$, respectively,  have $Rw^*wc-$property and $Lw^*wc-$property. If $A^{**}$ is weakly amenable, then $A$ is weakly amenable.
\begin{proof} Assume that $a^{\prime\prime}\in A^{**}$ and $(a_\alpha)_\alpha \subseteq A$ such that $a_\alpha \stackrel{w^*} {\rightarrow}a^{\prime\prime}$. Then for each $a^\prime\in A^*$, we have $a_\alpha a^\prime\stackrel{w^*} {\rightarrow}a^{\prime\prime}a^\prime$ in $A^*$. Since $A^*$ has $Rw^*wc-$property, $a_\alpha a^\prime\stackrel{w} {\rightarrow}a^{\prime\prime}a^\prime$ in $A^*$. Then for every $x^{\prime\prime}\in A^{**}$, we have
$$\langle   x^{\prime\prime}a_\alpha , a^\prime \rangle=\langle   x^{\prime\prime},a_\alpha  a^\prime \rangle\rightarrow
\langle   x^{\prime\prime},a^{\prime\prime}  a^\prime \rangle=\langle   x^{\prime\prime}a^{\prime\prime},  a^\prime \rangle.$$
It follows that $x^{\prime\prime}a_\alpha \stackrel{w^*} {\rightarrow}x^{\prime\prime}a^{\prime\prime}$. Since $A^{**}$ has $Lw^*wc-$property with respect to $A$, $x^{\prime\prime}a_\alpha \stackrel{w} {\rightarrow}x^{\prime\prime}a^{\prime\prime}$. If $D:A\rightarrow A^*$ is a bounded derivation, we extend it to a bounded linear mapping  $D^{\prime\prime}$ from $A^{**}$ into $A^{***}$. Suppose that  $a^{\prime\prime}, b^{\prime\prime}\in A^{**}$ and $(a_\alpha)_\alpha , (b_\beta)_\beta\subseteq A$ such that $a_\alpha \stackrel{w^*} {\rightarrow}a^{\prime\prime}$ and $b_\beta \stackrel{w^*} {\rightarrow}b^{\prime\prime}$. Since  $x^{\prime\prime}a_\alpha \stackrel{w} {\rightarrow}x^{\prime\prime}a^{\prime\prime}$ for every $x^{\prime\prime}\in A^{**}$, we have
$$\lim_\alpha \langle   D^{\prime\prime}(b^{\prime\prime}), x^{\prime\prime}a_\alpha\rangle=\langle   D^{\prime\prime}(b^{\prime\prime}), x^{\prime\prime}a^{\prime\prime}\rangle.$$
In the following we take limit on the $weak^*$ topologies.  Thus we have
$$\lim_\alpha \lim_\beta D(a_\alpha)b_\beta=D^{\prime\prime}(a^{\prime\prime})b^{\prime\prime}.$$
Consequently, we have
$$D^{\prime\prime}(a^{\prime\prime}b^{\prime\prime})=\lim_\alpha \lim_\beta D(a_\alpha b_\beta)=\lim_\alpha \lim_\beta D(a_\alpha)b_\beta +\lim_\alpha \lim_\beta a_\alpha D(b_\beta)$$$$=D^{\prime\prime}(a^{\prime\prime})b^{\prime\prime}+ a^{\prime\prime}D^{\prime\prime}(b^{\prime\prime}).$$
Since $A^{**}$ is weakly amenable, there is $a^{\prime\prime\prime}\in A^{***}$ such that $D^{\prime\prime}=\delta_{a^{\prime\prime\prime}}$. We conclude that
$D=D^{\prime\prime}\mid_{A}=\delta_{a^{\prime\prime\prime}}\mid_{A}$. Hence for each $x^\prime\in A ^*$, we have $D=x^\prime a^{\prime\prime\prime}\mid_{A }-a^{\prime\prime\prime}\mid_{A }x^\prime$. Take $a^\prime=a^{\prime\prime\prime}\mid_{A }$. It follows that $H^1(A,A^*)=0$.\end{proof}

\vspace{0.4cm}

\noindent{\it{\bf Theorem 3-5.}} Let $A$ be a Banach algebra and suppose that $D:A\rightarrow A^*$ is a surjective  derivation. If $D^{\prime\prime}$ is a derivation, then we have the following assertions.
\begin{enumerate}
\item ~$A^*$ and $A^{**}$, respectively,  have $w^*wc-$property and $Lw^*wc-$property with respect to $A$.
\item ~For every $a^{\prime\prime}\in A^{**}$, the mapping $x^{\prime\prime}\rightarrow a^{\prime\prime}x^{\prime\prime}$ from $A^{**}$ into $A^{**}$ is $weak^*-weak$ continuous.
\item ~ $A$ is Arens regular.
\item ~ If $A$ has $LBAI$, then $A$ is reflexive.\\
 \end{enumerate}

\begin{proof}
\begin{enumerate}
\item ~Since $D$ is surjective, $D^{\prime\prime}$ is surjective, and so by using [19, Theorem 2.2], we have $A^{***}A^{**}\subseteq D^{\prime\prime}(A^{**})A^{**}\subseteq A^*$. Suppose that $a^{\prime\prime}\in A^{**}$ and  $(a_\alpha)_\alpha \subseteq A$ such that
    $a_\alpha \stackrel{w^*} {\rightarrow}a^{\prime\prime}$. Then for each   $x^\prime\in A^*$, we have  $x^\prime a_\alpha \stackrel{w^*} {\rightarrow}x^\prime a^{\prime\prime}$. Since $A^{***}A^{**}\subseteq A^*$, $x^\prime a^{\prime\prime}\in A^*$. Then for every  $x^{\prime\prime}\in A^{**}$, we have
$$\langle   x^{\prime\prime}, x^\prime a_\alpha\rangle=\langle   x^{\prime\prime} x^\prime , a_\alpha\rangle\rightarrow   \langle   a^{\prime\prime} , x^{\prime\prime}x^\prime \rangle=\langle   x^\prime a^{\prime\prime} , x^{\prime\prime} \rangle=\langle    x^{\prime\prime},x^\prime a^{\prime\prime}  \rangle.$$
 It follows that  $x^\prime a_\alpha \stackrel{w} {\rightarrow}x^\prime a^{\prime\prime}$ in $A^*$. Thus  $x^\prime$ has  $Lw^*wc-$property with respect to $A$. The proof that   $x^\prime$ has  $Rw^*wc-$property with respect to $A$ is similar, and so $A^*$ has $w^*wc-$property. \\
Suppose that $x ^{\prime\prime\prime} \in A^{***}$.  Since $A^{***}A^{**}\subseteq A^*$,   $x^{\prime\prime}a_\alpha \stackrel{w^*} {\rightarrow}x^{\prime\prime}a^{\prime\prime}$ for each $x^{\prime\prime}\in A^{**}$. Then
$$\langle   x ^{\prime\prime\prime}, x^{\prime\prime}a_\alpha\rangle=\langle   x ^{\prime\prime\prime} x^{\prime\prime}, a_\alpha\rangle\rightarrow \langle   x ^{\prime\prime\prime} x^{\prime\prime}, a^{\prime\prime}\rangle=\langle   x ^{\prime\prime\prime}, x^{\prime\prime} a^{\prime\prime}\rangle.$$
It follows that $x^{\prime\prime}a_\alpha \stackrel{w} {\rightarrow}x^{\prime\prime}a^{\prime\prime}$.  Thus  $x^{\prime\prime}$ has  $Lw^*wc-$property with respect to $A$.
\item Suppose that $(a^{\prime\prime}_\alpha)_\alpha \subseteq A^{**}$ and  $a^{\prime\prime}_\alpha \stackrel{w^*} {\rightarrow}a^{\prime\prime}$. Let $x^{\prime\prime}\in A^{**}$. Then for every $x ^{\prime\prime\prime} \in A^{***}$, since $A^{***}A^{**}\subseteq A^*$, we have
$$\langle   x ^{\prime\prime\prime}, x^{\prime\prime}a^{\prime\prime}_\alpha\rangle=\langle   x ^{\prime\prime\prime} x^{\prime\prime}, a^{\prime\prime}_\alpha\rangle\rightarrow \langle   x ^{\prime\prime\prime} x^{\prime\prime}, a^{\prime\prime}\rangle=
\langle   x ^{\prime\prime\prime}, x^{\prime\prime} a^{\prime\prime}\rangle.$$

\item ~~It  follows from (2).
\item ~Let $(e_\alpha)_\alpha \subseteq A$ be a $BLAI$ for $A$. Then without loss generality, let $e^{\prime\prime}$ be a left unit for $A^{**}$ such that  $e_\alpha \stackrel{w^*} {\rightarrow}e^{\prime\prime}$. Suppose that $(a^{\prime\prime}_\alpha)_\alpha \subseteq A^{**}$ and  $a^{\prime\prime}_\alpha \stackrel{w^*} {\rightarrow}a^{\prime\prime}$. Then for every $a ^{\prime\prime\prime} \in A^{***}$, since $A^{***}A^{**}\subseteq A^*$, we have
$$\langle   a ^{\prime\prime\prime},a^{\prime\prime}_\alpha\rangle= \langle   a ^{\prime\prime\prime},e^{\prime\prime}a^{\prime\prime}_\alpha\rangle=\langle   a ^{\prime\prime\prime}e^{\prime\prime}, a^{\prime\prime}_\alpha\rangle\rightarrow \langle   a ^{\prime\prime\prime}e^{\prime\prime}, a^{\prime\prime}\rangle=\langle   a ^{\prime\prime\prime}, a^{\prime\prime}\rangle.$$
 It follows that   $a^{\prime\prime}_\alpha \stackrel{w} {\rightarrow}a^{\prime\prime}$. Consequently $A$ is reflexive.\end{enumerate}\end{proof}

\vspace{0.4cm}

\noindent{\it{\bf Corollary 3-6.}} Let $A$ be a Banach algebra and suppose that $D:A\rightarrow A^*$ is a surjective  derivation. Then the  following statements are equivalent.
\begin{enumerate}
\item ~$A^*$ and $A^{**}$, respectively,  have $Rw^*wc-$property and $Lw^*wc-$property.
\item ~For every $a^{\prime\prime}\in A^{**}$, the mapping $x^{\prime\prime}\rightarrow a^{\prime\prime}x^{\prime\prime}$ from $A^{**}$ into $A^{**}$ is $weak^*-weak$ continuous.\\\\
 \end{enumerate}

\noindent{\it{\bf Problems.}}\\ 1. Let $G$ be a locally compact group. What can  say  for the following sets?
\begin{enumerate}
\item $Z^\ell_{L^1(G)^{**}}((L^1(G)\hat{\otimes}L^1(G))^{**})=?~~~ , ~~~Z^r_{L^1(G)^{**}}((L^1(G)\hat{\otimes}L^1(G))^{**})=?$
\item  $Z^\ell_{(L^1(G)\hat{\otimes}L^1(G))^{**}}(L^1(G)^{**})=?~~~~ ~~, ~~~~Z^r_{(L^1(G)\hat{\otimes}L^1(G))^{**}}(L^1(G)^{**})=?$
\item  $Z^\ell_{L^1(G)^{**}}(L^1(G)^{**}\hat{\otimes}L^1(G)^{**})=?~~~~ , ~~~~Z^r_{L^1(G)^{**}}(L^1(G)^{**}\hat{\otimes}L^1(G)^{**})=?$

\end{enumerate}
2.  Suppose that  $S$ is a compact semigroup. Dose $L^1(S)^*$ and $M(S)^*$ have   $Lw^*wc-$property or $Rw^*wc-$property?\\\\

\bibliographystyle{amsplain}

\end{document}